\documentclass[12pt,reqno]{amsart}
\usepackage[margin=1in]{geometry}
\usepackage{amsmath,amssymb,amsthm,graphicx,amsxtra, setspace}
\usepackage[utf8]{inputenc}
\usepackage{mathrsfs}
\usepackage{hyperref}
\usepackage{xcolor}
\usepackage{upgreek}
\usepackage{mathtools}
\usepackage[numbers]{natbib}
\allowdisplaybreaks

\newtheorem{theorem}{Theorem}[section]
\newtheorem{lem}{Lemma}[section]

\newtheorem{Def}{Definition}[section]

\newtheorem{Ex}{Example}[section]
\newtheorem{Ass}{Assumption}[section]
\usepackage{latexsym}
\usepackage{hyperref}
\usepackage{graphicx}
\usepackage{epstopdf}

\allowdisplaybreaks

\let\originalleft\left
\let\originalright\right
\renewcommand{\left}{\mathopen{}\mathclose\bgroup\originalleft}
\renewcommand{\right}{\aftergroup\egroup\originalright}

\newcommand{\Addresses}{{
		\footnote{

			\noindent	 \textsuperscript{1,2} Department of Applied Mathematics and Scientific Computing, Indian Institute of Technology Roorkee, Roorkee, 247667, India.	
			
			\noindent  \textit{e-mail\textsuperscript{1}:} \texttt{g\_gupta@as.iitr.ac.in}
			
			\noindent  \textit{e-mail\textsuperscript{2}:} \texttt{jay.dabas@gmail.com.}

			\noindent \textsuperscript{*}Corresponding author.
			
			\textit{Key words:}  Hemivariational inequalities; Approximate controllability; Riemann-Liouville fractional derivative; multivalued analysis; Strongly continuous semigroup.
			
			Mathematics Subject Classification (2020): 34G25, 34A08, 49J40 , 37L05, 93B05. 
			
}}}

\begin{document}
	\title[Controllability of Fractional Order Hemivariational Inequalities]{Existence And Approximate Controllability for a class of Fractional Order Hemivariational Inequalities.\Addresses}
	\author [Garima Gupta AND Jaydev Dabas]{Garima Gupta\textsuperscript{1} AND  Jaydev Dabas\textsuperscript{2*}}
	\maketitle
	
	\begin{abstract}
		This paper discusses the approximate controllability of a fractional differential control problem driven by a nonlinear hemivariational inequality in a Hilbert space. First, we prove the existence of a mild solution for a fractional control inclusion problem which is equivalent to a hemivariational inequality by using the nonsmooth analysis and fixed point technique. Further, we established sufficient conditions for the approximate controllability of our inclusion problem by taking corresponding linear system is approximately controllable. The existence and controllability results obtained for the inclusion problem are valid for considered nonlinear hemivariational problem. Finally, we provide an example to illustrate the efficiency of the developed results.
	\end{abstract}
	\section{\textbf{Introduction}}\label{intro}\setcounter{equation}{0}
	\vspace{0.5cm}
	
	Hemivariational inequalities are a generalization of variational inequalities that arise in the study of nonconvex and nonsmooth energy functions. They have various applications across different fields due to their ability to model complex phenomena involving nonconvex and nonsmooth potentials.
	Panagiotopoulos initially introduced the concept of hemivariational inequality in 1981(cf. \cite{PANAGIOTOPOULOS1981335}). He utilized hemivariational inequalities to address mechanical problems characterized by nonconvex and nonsmooth superpotentials, see for example,
	\cite{panagiotopoulos1993hemivariational,panagiotopoulos2012inequality}. Over time, an increasing number of scholars have made significant contributions to the exploration of solution existence in hemivariational inequalities and various authors have proven the existence of solutions for hemivariational inequalities under different assumptions and hypotheses. For detailed information of existence of solution and its nature we refer Kavitha, Vijayakumar, Shukla, Nisar, Kottakkaran and Udhayakumar\cite{kavitha2021results}, Mohan Raja, Vijayakumar, Udhayakumar, Nisar and Kottakkaran \cite{mohan2024results}, Jiang, Wei, Zhouchao, Guoji and Irene\cite{jiang2023topological}, Zeng, Liu and Migorski\cite{zeng2018class}, Ma, Dineshkumar, Vijayakumar, Udhayakumar, Shukla, Anurag and Kottakkaran \cite{ma2023hilfer} and the references given in these articles.
	The notion of noninteger derivatives and integrals represents an extension of the conventional calculus based on integer orders. This extension is motivated by the distinctive memory-like characteristics inherent in fractional derivatives, rendering them more suitable for describing the properties of diverse real materials compared to their integer-order counterparts. Over the past two decades, fractional calculus has drawn the interest of physicists, mathematicians, and engineers, leading to notable contributions in both theoretical advancements and practical applications of fractional differential equations. For more comprehensive insights into fractional calculus and fractional differential equations, readers are directed to the monograph authored by Kilbas\cite{kilbas2006theory}. Hemivariational inequalities with fractional derivatives are essential in modeling anomalous diffusion processes where the standard diffusion equations fail, such as in porous media or heterogeneous materials. The specifications of initial conditions for Riemann–Liouville fractional derivatives or integrals are pivotal in addressing certain practical challenges. Heymans and Podlubny\cite{heymans2006physical,podlubny2001geometric} have illustrated that it is feasible to assign a physical significance to initial conditions formulated using Riemann–Liouville fractional derivatives or integrals, particularly in the realm of viscoelasticity. Such initial conditions are deemed more suitable than those that are physically interpretable.
	
	The nonlocal initial condition proves to be more effective in physics compared to the classical initial condition $u(0) = u_0$. To illustrate, in 1993, Deng \cite{deng1993exponential} utilized the nonlocal condition to characterize the diffusion phenomenon of a small amount of gas within a transparent tube. In this context, condition (1.2) facilitates additional measurements at $t_k$, where $k = 1, 2, \ldots, m$, offering greater precision than measurements solely at $t = 0$. Furthermore, in 1999, Byszewski\cite{byszewski1999existence} highlighted that if $c_k \neq 0$, where $k = 1, 2, \ldots, m$, the outcomes can be employed in kinematics to ascertain the evolutionary path $t \rightarrow u(t)$ of a physical object. This is particularly useful when the positions $u(0), u(t_1), \ldots, u(t_m)$ are unknown, but the nonlocal condition (1.2) is confirmed to hold. Few more articles by Mahmudov\cite{mahmudov2008approximate}, Wang\cite{wang2017approximate} and Chen\cite{chen2014existence,chen2020existence} considered semilinear systems with non-local conditions and proved the exixtence of solution.
	
	The introduction of controllability by Kalman in 1963\cite{kalman1963mathematical} marked the inception of an active research area, owing to its significant applications in physics. There are various works on approximate controllability of systems represented by fractional differential equations, integrodifferential equations, differential inclusions,
	neutral functional differential equations, and impulsive differential equations in Banach spaces; see\cite{arora2020approximate,liu2013controllability,mahmudov2008approximate} and their references.
	In recent years, the exploration of control systems governed by Caputo fractional evolution equations has seen considerable attention (see\cite{chen2020existence,wang2017approximate,mohan2024results,shi2016study,zhou2013existence}). Despite this, the topic of approximate controllability for fractional evolution differential equations with Riemann–Liouville fractional derivative with local and nonlocal initial conditions under different hypothesis has been studied by many authors. For reference see the literature\cite{du2011initialized,liang2022existence,liu2013approximate,liu2013controllability,liu2015approximate,shu2019approximate}. This gap in knowledge serves as the motivation for the present work. The objective of this paper is to present suitable sufficient conditions for the existence and  approximate controllability of fractional differential Hemivariational inequalities involving Riemann–Liouville fractional derivatives.
	
	Let $\mathbb{H}$  be a seperable Hilbert space and $\mathbb{U}$ be a Hilbert space. In this work, we investigate the existence of a mild solution and the approximate controllability of the following semilinear fractional differential hemivariational inequality:
	\begin{equation}\label{P}
		\begin{cases}
			\langle-^R\mathrm{D_{0,t}^{\alpha}} x(t)+\mathrm{A}x(t)+\mathrm{B}u(t),v\rangle_{\mathbb{H}}+F^{0}(t,x(t);v)\ge0,\quad t\in J=[0,b],\forall v\in \mathbb{H},\\
			{I_{0,t}}^{1-\alpha}x(t)|_{t=0}=\sum _{k=1}^{m}c_{k}x_{k}.
		\end{cases}
	\end{equation}
	where,
	$\langle.,.\rangle_{\mathbb{H}}$ denotes the scalar product of the separable Hilbert space $\mathbb{H}$ and the norm in $\mathbb{H}$ is denoted by $\|.\|_{\mathbb{H}}$ , $^R\mathrm{D_{0,t}^{\alpha}}$ denotes the Riemann-Liouville fractional derivative of order $\alpha \in(0,1)$ with the lower limit zero and ${I_{0,t}^{1-\alpha}}$ denotes the Riemann-Liouville fractional integral of order $1-\alpha$ with lower limit zero, $\mathrm{A}:D(\mathrm{A})\subseteq\mathbb{H}\rightarrow\mathbb{H}$ is the infinitesimal generator of a $C_{0}-$semigroup $\mathcal{T}(t)(t\ge0)$ on $\mathbb{H}$. For $\alpha>\frac{1}{2}$ the control function $u$ takes value in $L^{2}(J,\mathbb{U})$ of admissible control functions for a hilbert space $\mathbb{U}$ , $\mathrm{B}:\mathbb{U}\rightarrow\mathbb{H}$ is a bounded linear operator,
	$ F^{0}(t,.;.)$ stands for the generalized Clarke directional derivative of a locally Lipschitz function $F(t,.):\mathbb{H}\rightarrow\mathbb{R}$,
	$0<t_{1}<t_{2}<......<t_{m}<b, m\in\mathbb{N} $, $c_{k}$ are real constant, $c_{k}\ne0, k=1,2,....,m$ and $x_{k}=x(t_{k}) $ for $k=1,2,....,m$.
	
	\section{\textbf{Preliminaries}}\label{pre}\setcounter{equation}{0}
	In this section, we recall some fundamental definitions, notations, which will help us to establish existence and controllability result for the system \eqref{P}.
	
	The norm of a Banach space $\mathbb{X}$ will be denoted by $\|\cdot\|_{\mathbb{X}} \cdot L_{b}(\mathbb{X}, \mathbb{Y})$ denotes the space of bounded linear operators from a Banach space $\mathbb{X}$ to Banach space $\mathbb{Y}$. For the uniformly bounded $C_{0}$-semigroup $\mathcal{T}(t)(t \geq 0)$, we set $M:=\sup _{t\ge0}\|T(t)\|_{L_{b}(\mathbb{X})}<\infty$. Let $C(J, X)$ denote the Banach space of all $\mathbb{X}$ value continuous functions from $J=[0, b]$ to $\mathbb{X}$ with the norm $\|x\|_{C}=\sup _{t \in J}\|x(t)\|_{\mathbb{X}}$. Let $C_{1-\alpha}(J, \mathbb{X})=\left\{x: t^{1-\alpha} x(t) \in C(J, \mathbb{X})\right\}$ with the norm
	
	$$
	\|x\|_{C_{1-\alpha}}=\sup \left\{t^{1-\alpha}\|x(t)\|_{\mathbb{X}}: t \in J\right\} .
	$$
	
	Obviously, the space $C_{1-\alpha}(J, \mathbb{X})$ is a Banach space.
	Some definitions related to fractional integral and derivatives are as follows:
	\begin{Def}\cite{kilbas2006theory}\label{def2.1}
		The fractional integral of a function $z:[a,b]\to\mathbb{R}$, $a,b\in\mathbb{R}$ with $a<b$, of order $\alpha>0$ is defined as
		\begin{align*}
			\mathrm{I}_{a,t}^{\alpha}z(t):=\frac{1}{\Gamma(\alpha)}\int_{a}^{t}\frac{z(s)}{(t-s)^{1-\alpha}}\mathrm{d}s,\ \mbox{ for a.e. } \  t\in[a,b],
		\end{align*}
		
		where $z\in\mathrm{L}^1([a,b];\mathbb{R})$ and $\Gamma(\alpha)=\int_{0}^{\infty}t^{\alpha-1}e^{-t}\mathrm{d}t$ is the Euler gamma function.
	\end{Def}
	\begin{Def}\cite{kilbas2006theory}\label{def2.2}
		The Riemann-Liouville fractional derivative of a function $z:[a,b]\to\mathbb{R}$  of order $\alpha>0$ is given as 
		\begin{align*}
			^R\mathrm{D}_{a,t}^{\alpha}z(t):=\frac{1}{\Gamma(n-\alpha)}\frac{d^n}{dt^n}\int_{a}^{t}(t-s)^{n-\alpha-1}z(s)\mathrm{d}s,\ \mbox{ for a.e. }\ t\in[a,b],
		\end{align*}
		with $n-1< \alpha<n$.
		
	\end{Def}	
	
	If $z$ is an abstract function with values in $\mathbb{X}$, then the integrals which appear in \ref{def2.1} and \ref{def2.2} are taken in Bochner's sense, that is: a measurable function $z$ maps from $\left[0,+\infty\right)$  to $\mathbb{X}$ is Bochner integrable if $\|z\|$ is Lebesgue integrable.
	Furthermore, given a Banach space $\mathbb{X}$, we will use the following notations.
	
	\begin{align*}
		& \mathcal{P}_{cl,cv}(X):=\{\Omega \subseteq \mathbb{X}: \Omega \ \text {is nonempty, closed (convex)}\}, \\
		& \mathcal{P}_{(w)cp(cv)}(X):=\{\Omega \subseteq \mathbb{X}: \Omega \ \text {is nonempty, (weakly) compact (convex)}\} .
	\end{align*}
	
	Now, we introduce some basic definitions and results from multivalued analysis. For more details on multivalued maps, please see the book\cite{deimling2011multivalued}.
	\begin{enumerate}
		\item[\textit{(i)}]For a given Banach space $\mathbb{X}$, a multivalued map $F: \mathbb{X} \rightarrow 2^{\mathbb{X}} \backslash\{\emptyset\}:=\mathcal{P}(\mathbb{X})$ is convex (closed) valued, if $F(x)$ is convex (closed) for all $x \in \mathbb{X}$.
		
		\item[\textit{(ii)}] $F$ is called upper semicontinuous (u.s.c. for short) on $\mathbb{X}$, if for each $x \in \mathbb{X}$, the set $F(x)$ is a nonempty, closed subset of $\mathbb{X}$, and if for each open set $V$ of $\mathbb{X}$ containing $F(x)$, there exists an open neighborhood $N$ of $x$ such that $F(N) \subseteq V$.
		
		\item[\textit{(iii)}] $F$ is said to be completely continuous if $F(V)$ is relatively compact, for every bounded subset $V \subseteq \mathbb{X}$.
		
		\item[\textit{(iv)}] Let $\Sigma$ is the $\sigma$- algebra of subsets of the set $\Omega$, $(\Omega, \Sigma)$ be a measurable space and $(\mathbb{X}, d)$ a separable metric space. A multivalued map $F: \Omega \rightarrow \mathcal{P}(\mathbb{X})$ is said to be measurable, if for every closed set $C \subseteq \mathbb{X}$, we have $F^{-1}(C)=\{t \in \Omega: F(t) \cap C \neq \emptyset\} \in \Sigma$.	
	\end{enumerate}
	Now we recall the few elements of nonsmooth analysis(see \cite{clarke1983optimization} for detailed information).
	\begin{Def}
		Let $h:\mathbb{X}\longrightarrow\mathbb{R}$ be a locally Lipschitx function on a Banach space $\mathbb{X}$. The generalized directional derivative of $h$ at $y\in \mathbb{X}$ in the direction $z\in \mathbb{X}$ is defined by 
		\begin{equation*}
			h^{0}(y;z):=\lim_{\lambda\rightarrow0^{+}}\sup_{\eta\rightarrow y}\frac{h(\eta+\lambda z)-h(\eta)}{\lambda}.
		\end{equation*}
	\end{Def}
	The generalized gradient of $h$ at $y\in X$ is the subset of $\mathbb{X}^{*}$ which is the dual space of $\mathbb{X}$, is given by $$ \partial h(y):=\{y^{*}\in \mathbb{X}^{*}: h^{0}(y;z)\ge \langle y^{*},z\rangle \forall z\in \mathbb{X}\},$$
	
	Now we Consider the following semilinear inclusion
	\begin{equation}\label{P2}
		\begin{cases}
			^R\mathrm{D_{0,t}^{\alpha}}x(t)\in \mathrm{A}x(t)+\mathrm{B}u(t)+\partial F(t,x(t)), \quad t\in J=[0,b] ,\\
			{I_{0,t}^{1-\alpha}}x(t)|_{t=0}=\sum _{k=1}^{m}c_{k}x_{k},
		\end{cases}	
	\end{equation}
	where, $\partial F$ is the generalized Clarke subdifferential of a locally Lipschitz function $F(t,.):\mathbb{H}\rightarrow\mathbb{R}$. If $x\in C_{1-\alpha}(J,\mathbb{H})$ is a solution of \eqref{P2}, then there exists $f(t)\in \partial F(t,x(t))$ such that $f(t)\in L^{1}(J,\mathbb{H})$ and
	\begin{equation*}
		\begin{cases}
			^R\mathrm{D_{0,t}^{\alpha}}x(t)= \mathrm{A}x(t)+\mathrm{B}u(t)+f(t),\quad t\in J=[0,b] ,\\
			{I_{0,t}^{1-\alpha}}x(t)|_{t=0}=\sum _{k=1}^{m}c_{k}x_{k},
		\end{cases}	
	\end{equation*}
	which implies
	\begin{equation*}
		\begin{cases}
			\langle-^R\mathrm{D_{0,t}^{\alpha}} x(t)+\mathrm{A}x(t)+\mathrm{B}u(t),v\rangle_{\mathbb{H}}+\langle f(t),v\rangle_{\mathbb{H}}=0,\quad t\in J=[0,b], \forall v\in \mathbb{H}\\
			{I_{0,t}^{1-\alpha}}x(t)|_{t=0}=\sum _{k=1}^{m}c_{k}x_{k}.
		\end{cases}
	\end{equation*}
	since $f\in \partial F(t,x(t))$ and $\langle f(t),v\rangle_{\mathbb{H}}\le F^{0}(t,x(t);v)$, we obtain
	\begin{equation*}
		\begin{cases}
			\langle-^R\mathrm{D_{0,t}^{\alpha}} x(t)+\mathrm{A}x(t)+\mathrm{B}u(t),v\rangle_{\mathbb{H}}+F^{0}(t,x(t);v)\ge0,\quad t\in J=[0,b], \forall v\in \mathbb{H}\\
			{I_{0,t}^{1-\alpha}}x(t)|_{t=0}=\sum _{k=1}^{m}c_{k}x_{k}.
		\end{cases}
	\end{equation*}
	Therefore, in order to study the hemivariational inequality \eqref{P}, we only need to deal with the semilinear inclusion \eqref{P2}.\newline
	Further, we define the operator
	\begin{align*}
		\mathcal{T_\alpha}(t) = \alpha\int_{0}^{\infty}\theta\xi_{\alpha}(\theta)\mathcal{T}(t^{\alpha}\theta)\mathrm{d}\theta,
	\end{align*}	
	\begin{align*}
		\xi_{\alpha}(\theta) = \frac{1}{\alpha}\theta^{-1-(1/\alpha)}\omega_{\alpha}(\theta^{-1/\alpha)},
	\end{align*}
	\begin{align*}
		\omega_{\alpha}(\theta) = \frac{1}{\pi}	\sum_{n = 1}^{\infty}(-1)^{n-1}\theta^{-n\alpha-1}\frac{\Gamma(n\alpha+1)}{n!}sin(\pi n\alpha), \quad \theta \in (0, \infty).
	\end{align*}
	\begin{Ass}\label{ass1}
		$\sum_{k=1}^{m}\left| c_{k}t_{k}^{\alpha-1}\right|<\frac{\Gamma{\alpha}}{M}$.
	\end{Ass}
	
	From assumption \ref{ass1} , we have
	\begin{align}\label{eqn2.2}
		\left\|\sum_{k=1}^{m}c_{k}t_{k}^{\alpha-1}\mathcal{T_\alpha}(t_{k})\right\|<1.
	\end{align} 
	By equation\eqref{eqn2.2} and operator spectrum theorem, we know that 
	\begin{align}\label{eqn2.3}
		\mathcal{O}=\big(I-\sum_{k=1}^{m}c_{k}t_{k}^{\alpha-1}\mathcal{T_\alpha}(t_{k})\big)^{-1},
	\end{align}
	exists and is a bounded operator with $D(\mathcal{O})=\mathbb{H}$. Furthermore, by Neumann expression, $\mathcal{O}$ can be expressed by 
	\begin{align*}
		\mathcal{O}=\sum_{n=0}^{\infty}\big(\sum_{k=1}^{m}c_{k}t_{k}^{\alpha-1}\mathcal{T_\alpha}(t_{k})\big)^{n}.
	\end{align*}
	Therefore,
	\begin{align*}
		\left\|\mathcal{O}\right\|&\le\sum_{n=0}^{\infty}\left\|\sum_{k=1}^{m}c_{k}t_{k}^{\alpha-1}\mathcal{T_\alpha}(t_{k})\right\|^{n},\\
		&\le \frac{1}{1-\frac{M}{\Gamma\alpha}\sum_{k=1}^{m}\left| c_{k}t_{k}^{\alpha-1}\right|}.
	\end{align*}
	By the above discussion, \cite{liu2015approximate} and \cite{liu2013approximate} We know that the mild solution for the fractional inclusion problem\eqref{P2} can be written as
	\begin{align}\label{eqn2.4}
		x(t)=t^{\alpha-1}\mathcal{T_\alpha}(t)\left({I_{0,t}}^{1-\alpha}x(t)|_{t=0}\right)+\int_{0}^{t}(t-s)^{\alpha-1}\mathcal{T_\alpha}(t-s)[\mathrm{B}u(s)+f(s)]\mathrm{d}s,
	\end{align}
	From \eqref{eqn2.4} we have for each $t_{k}$
	\begin{align}\label{eqn2.5}
		x(t_{k})= {t_{k}}^{\alpha-1}\mathcal{T_\alpha}(t_{k})\left({I_{0,t}}^{1-\alpha}x(t)|_{t=0}\right)+\int_{0}^{t_{k}}(t_{k}-s)^{\alpha-1}\mathcal{T_\alpha}(t_{k}-s)[\mathrm{B}u(s)+f(s)]\mathrm{d}s.
	\end{align}
	Using Assumption \ref{ass1} and the estimates\eqref{eqn2.3}, \eqref{eqn2.4} and \eqref{eqn2.5},  we get
	\begin{align}\label{eqn2.6}
		{I_{0,t}}^{1-\alpha}x(t)|_{t=0}=\sum_{k=1}^{m}c_{k} \mathcal{O}\int_{0}^{t_{k}}{(t_{k}-s)}^{\alpha-1} \mathcal{T_\alpha}(t_{k}-s)[\mathrm{B}u(s)+f(s)]\mathrm{d}s.
	\end{align}
	By \eqref{eqn2.4} and \eqref{eqn2.6}, we can write 
	\begin{align}\label{eqn2.7}
		x(t)&=\sum_{k=1}^{m}c_{k} t^{\alpha -1}\mathcal{T_\alpha} \mathcal{O}\int_{0}^{t_{k}}{(t_{k}-s)}^{\alpha-1} \mathcal{T_\alpha}(t_{k}-s)[\mathrm{B}u(s)+f(s)]\mathrm{d}s
		\nonumber\\&\quad+\int_{0}^{t}{(t-s)}^{\alpha-1} \mathcal{T_\alpha}(t-s)[\mathrm{B}u(s)+f(s)]\mathrm{d}s.
	\end{align}
	For convenience, we introduce the function $G(t,s)$ as follows:
	\begin{align}\label{eqn2.8}
		G(t,s)=\sum_{k=1}^{m}\chi_{t_k}t^{\alpha-1}\mathcal{T_\alpha}(t)(t_{k}-s)^{\alpha-1}\mathcal{O}\mathcal{T_\alpha}(t_{k}-s)+\chi_{t}(s)(t-s)^{\alpha-1}\mathcal{T_\alpha}(t-s),
	\end{align}
	with	
	\begin{equation*}
		\chi_{t_k}(s)=
		\begin{cases}
			c_k,\hspace{.2cm} s\in\left[0,t_{k}\right) \\ 0,\hspace{.2cm} s\in[t_{k},b],
		\end{cases}
	\end{equation*}	
	\begin{equation*}
		\chi_{t}(s)=
		\begin{cases}
			1,\hspace{.2cm} s\in\left[0,t\right) \\ 0,\hspace{.2cm} s\in[t,b].
		\end{cases}
	\end{equation*}
	Therefore, by \eqref{eqn2.7} and \eqref{eqn2.8} we know that the solution of fractional inclusion\eqref{P2} can also be expressed as
	\begin{align}
		x(t)=\int_{0}^{b}G(t,s)[\mathrm{B}u(s)+f(s)]\mathrm{d}s.
	\end{align} 
	
	Now we may define a mild solution of problem \eqref{P2} as follows:
	\begin{Def}
		For each $u\in L^{2}(J,\mathbb{U})$, a function $x\in C_{1-\alpha}(J,\mathbb{H})$ is called a mild solution of the control system \eqref{P} if ${I_{t}^{1-\alpha}}x(t)|_{t=0}=\sum _{k=1}^{m}c_{k}x_{k}$ and there exists $f\in L^{1}(J,\mathbb{H})$ such that $f(t)\in\partial F(t,x(t))$ a.e. on $t\in J$ and
		\begin{equation}\label{MS}
			x(t)=\int_{0}^{b}G(t,s)[\mathrm{B}u(s)+f(s)]\mathrm{d}s.
		\end{equation}
	\end{Def}

	\begin{lem}\label{lem2.1}\cite{liu2015approximate}
		The operator $\mathcal{T_\alpha}(t)$ has the following properties:
		\begin{enumerate}
			\item For any fixed $t>0$, $\mathcal{T_\alpha}(t) $ is linear and bounded operator, that is for any $x\in \mathbb{H}$,
			\begin{equation}
				\left\|\mathcal{T_\alpha}(t)x\right\|\le \frac{M}{\Gamma\alpha}\left\|x\right\|,
			\end{equation}
			\item $\mathcal{T_\alpha}(t)(t\ge0)$ is strongly continuous.
		\end{enumerate}
	\end{lem}   
	\begin{Def}
		Let $x$ be a mild solution of system \eqref{P2} corresponding to the control $u\in L^{2}(J,\mathbb{U})$. Fractional evolution inclusion\eqref{P2} is said to be approximately controllable on the interval $J$ if the set $\overline{\mathcal{R}_{f}(b)}=\mathbb{H}$, where the set
		\begin{equation*}
			\mathcal{R}_{f}(b)=\{x(b)\in \mathbb{H}: u\in L^{2}(J,\mathbb{U}) \},
		\end{equation*}
		is called the reachable set of \eqref{P2}.
	\end{Def}  
	\section{\textbf{Existence of Mild Solution}}
	In this section we will prove the existence of mild solution of system \eqref{P2} by assuming some sufficient conditions and fixed point theorem. We start this section by defining the following operators
	\begin{align*}
		\Gamma_{0}^{b}=\int_{0}^{b}G(b,s) \mathrm{B} \mathrm{B}^{*}G^{*}(b,s)\mathrm{d}s,\quad \frac{1}{2}<\alpha \le 1,
	\end{align*}
	and
	\begin{align*}
		\mathrm{R}\left(a, \Gamma_{0}^{b}\right)=\left(a I+\Gamma_{0}^{b}\right)^{-1}, \quad a>0,
	\end{align*}
	where $B^{*}$,$\mathcal{O}^{*}$ and $\mathcal{T_\alpha}^{*}$ is the adjoint of $\mathrm{B}$, $\mathcal{O}$ and $\mathcal{T_\alpha}$ respectively, and $G^{*}$ is the adjoint of $G$ defines as:
	\begin{align*}
		G^{*}(b,s)=\sum_{k=1}^{m}\chi_{t_k}(s)t^{\alpha-1}\mathcal{T}_{\alpha}^{*}(t)\mathcal{O}^{*}(t_{k}-s)^{\alpha-1}\mathcal{T_{\alpha}}^{*}(t_{k}-s)+\chi_{t}(s)(t-s)^{\alpha-1}\mathcal{T_{\alpha}}^{*}(t-s).
	\end{align*}
	We consider the linear fractional differential control system:
	\begin{equation}\label{LP}
		\begin{cases}
			^R\mathrm{D_{0,t}^{\alpha}}x(t)= \mathrm{A}x(t)+\mathrm{B}u(t), \hspace{0.5cm}t\in J=[0,b] ,\quad \frac{1}{2}<\alpha \le 1,\\
			{I_{0,t}^{1-\alpha}}x(t)|_{t=0}=\sum _{k=1}^{m}c_{k}x_{k}.
		\end{cases}	
	\end{equation}
	\begin{lem}\label{lem3.1}\cite{mahmudov2003approximate}
		The linear fractional differential system \eqref{LP}  is approximately controllable on $J$ if and only if aR $\left(a, \Gamma_{0}^{b}\right) \rightarrow 0$ as $a \rightarrow 0^{+}$in the strong operator topology.
	\end{lem}
	\begin{lem}\label{lem3.2}\cite{kamenskii2011condensing}
		Let $\mathbb{X}$ be a Banach space. Let $F:J\times \mathbb{X}\longrightarrow\mathcal{P}_{cp,cv}(\mathbb{X})$ be an $L^{1}-$ Caratheodory multivalued map with $\mathcal{S}_{F(y)}=\left\{g\in L^{1}(J,\mathbb{X}): g(t)\in F(t,y(t)), \text{for a.e.}t\in J\right\}$ being nonempty and let $\Gamma $ be a linear continuous mapping from $L^{1}(J,\mathbb{X})$ to $C(J,\mathbb{X})$, then the operator .
		\begin{align*}
			\Gamma o \mathcal{S}_{F}:C(J,\mathbb{X})\longrightarrow\mathcal{P}_{cp,c}\left(C(J,\mathbb{X})\right),
		\end{align*}
		\begin{align*}
			y\longrightarrow(\Gamma o \mathcal{S}_{F})(y):=\Gamma(\mathcal{S}_{F(y)}),
		\end{align*}
		is a closed graph operator in $C(J,\mathbb{X})\times C(J,\mathbb{X})$.
	\end{lem}
	\begin{theorem}\label{thrm3.1}\cite{deimling2011multivalued}
		Let $D$ be a bounded, convex, and closed subset in the Banach space $\mathbb{X}$ and let $V: D \rightarrow 2^{X} \backslash\{\emptyset\}$ be a u.s.c. condensing multivalued map. If, for every $x \in D, V(x)$ is a closed and convex set in $D$, then $V$ has a fixed point.
	\end{theorem} 
	To prove the existence of mild solution we need the following assumptions:
	\begin{Ass}\label{ass3.1}
		\begin{enumerate}
			\item[(H1)] $\mathcal{T_\alpha}(t)$ is compact,
			\item[(H2)] the function $ t\mapsto F(t,x)$ is measurable for all $x \in \mathbb{H}$,
			\item[(H3)]the function $ x \mapsto F(t, x)$ is locally Lipschitz continuous for a.e. $t \in J$,
		\end{enumerate}
		
		\begin{enumerate}
			\item[(H4)] for each fixed $x\in C_{1-\alpha}(J,\mathbb{H})$ the set
			\begin{equation*}
				S_{\partial F,x} = \{f \in L^1(J, \mathbb{H}) : f(t) \in \partial F(t, x(t))\},
			\end{equation*}
			is nonempty,
			\item[(H5)] there exist a function $ P(t)\in L^{\frac{1}{\gamma}}(J, \mathbb{R}^+)$ with $0\le\gamma<\alpha$ and a nondecreasing function $\psi :\mathbb{R}\longrightarrow \mathbb{R}^{+}$, such that
			\begin{equation*}
				\|\partial F(t, x)\|_{\mathbb{H}} = \text{sup}\{\|f(t)\|_{\mathbb{H}} : f(t) \in \partial F (t,x)\} \le P(t)\psi(\|x\|_{D}),
			\end{equation*} 
			for any $t\in J$ for all $x\in \mathbb{H}$ and for each $r>0$, there exists $0<\rho<1$, such that
			\begin{align*}
				\lim _{r \rightarrow \infty} \inf \frac{\psi(r)}{r}\|P\|_{L^{2}}=\rho<1 .
			\end{align*}	
		\end{enumerate}
	\end{Ass}
	
	\begin{theorem}\label{Thrm3.2}
		If the assumption \ref*{ass1} and all the conditions (H1)-(H5) of assumption\ref{ass3.1} are satisfied, then the system \eqref{P2} has a mild solution. 
	\end{theorem}
	\textbf{Proof:} We consider a set
	\begin{align*}
		B_{r}=\left\{x \in C_{1-\alpha}(J, \mathbb{H}):\|x\| \leq r, r>0\right\}. 
	\end{align*}
	on the space $C_{1-\alpha}(J,\mathbb{H})$, We easily know that $B_{r}$ is a bounded, closed, and convex set in  $C_{1-\alpha}(J, \mathbb{H})$. For $a>0$, for all $x(\cdot) \in C_{1-\alpha}(J,\mathbb{H}), x_{1} \in \mathbb{H}$, we take the control function as
	\begin{align*}
		u(t)= \mathrm{B}^{*}G^{*}(b,t) \mathrm{R}\left(a, \Gamma_{0}^{b}\right) P(x(\cdot)),
	\end{align*}
	where
	\begin{align*}
		P(x(\cdot))=x_{1}-\int_{0}^{b}G(b,s) f(s)\mathrm{d}s, \quad f \in S_{\partial F, x} .
	\end{align*}
	By this control, we define the operator $\Phi_{a}: C_{1-\alpha}(J,\mathbb{H}) \rightarrow$ $\mathscr{P}\left(C_{1-\alpha}(J,\mathbb{H})\right)$ as follows:
	\begin{align*}
		\Phi_{a}(x)=\{\tau \in C_{1-\alpha}(J, \mathbb{H}): \tau(t)=\int_{0}^{b}G(t,s)[f(s)+\mathrm{B}u(s)]\mathrm{d}s, \quad f \in S_{\partial F, x}, t \in\left(0,b\right].
	\end{align*}
	To prove that the operator	$\Phi_{a}: C_{1-\alpha}(J,\mathbb{H}) \rightarrow$ $\mathscr{P}\left(C_{1-\alpha}(J,\mathbb{H})\right)$ has a fixed point, we subdivided the proof into following steps:
	
	\textbf{Step$1$:} $\Phi_{a}$ is convex for each $x\in C_{1-\alpha}(J,\mathbb{H})$.
	
	if $\tau_{1},\tau_{2}\in \Phi_{a}(x)$, then for each $t\in J$, $f_{1}$, $f_{2} \in S_{\partial F,x}$ s.t.
	\begin{align*}
		\tau_{i}(t)=\int_{0}^{b}G(t,s)[f_{i}(s)+\mathrm{B}\mathrm{B}^{*}G^{*}(b,t) \mathrm{R}\left(a, \Gamma_{0}^{b}\right)\{x_{1}-\int_{0}^{b}G(b,\mu)f_{i}(s) \,d\mu\}]\mathrm{d}s.
	\end{align*}
	Let $0\le \lambda \le1$, then for each $t\in J$, we have
	\begin{align*}
		&\lambda \tau_{1}(t)+(1-\lambda)\tau_{2}(t)\\&=\int_{0}^{b}G(t,s)[\lambda f_{1}(s)+(1-\lambda)f_{2}(s)]\mathrm{d}s+
		\int_{0}^{b}G(t,s)\mathrm{B}\mathrm{B}^{*}G^{*}(b,t) \mathrm{R}\left(a, \Gamma_{0}^{b}\right) \bigg[x_{1}\\&\qquad-\int_{0}^{b}G(b,\mu) \big(\lambda f_{1}(\mu)+(1-\lambda)f_{2}(\mu)\big) \,d\mu\bigg]\mathrm{d}s.
	\end{align*}
	Since $S_{\partial F,x}$ is convex $(\text{as} \partial F\text{has convex values})$, $\lambda f_{1}+(1-\lambda)f_{2}\in S_{\partial F,x}$, thus $\lambda \tau_{1}(t)+(1-\lambda)\tau_{2}(t)\in \Phi_{a}(x)$.
	%
	
	\textbf{Step$2$:} For each $a>0$, there is a positive constant $r_{0}=r(a)$, such that $\Phi_{a}(B_{r_{0}})\subseteqq B_{r_{0}}$. 
	
	If this is not true, there $\exists a>0$ such that $\forall r>0$ there exixts a $\stackrel{-}{x}$ such that $\Phi_{a}(\stackrel{-}{x})\nsubseteq$, that is
	\begin{equation*}
		\|\Phi_{a}(\stackrel{-}{x})\|=sup\{\|\tau\|_{C_{1-\alpha}(J,\mathbb{H})}: \tau \in \Phi_{a}(\stackrel{-}{x}) >r\}.
	\end{equation*}
	Since 
	\begin{align*}
		\tau(t)&=\int_{0}^{b} G(t,s)[f(s)+\mathrm{B}u(s)]\mathrm{d}s,\\
		\tau(t)&=\int_{0}^{b} G(t,s)f(s)\,ds+\int_{0}^{b}G(t,s)\mathrm{B}\mathrm{B}^{*}G^{*}(t,s) \mathrm{R}\left(a, \Gamma_{0}^{b}\right) \{x_{1}-\int_{0}^{b}G(b,\mu)f(\mu)d\mu\}\mathrm{d}s,
	\end{align*}
	for some $f\in S_{\partial F,\stackrel{-}{x}}$.
	
	By using Holder's inequality and (\textit{H5}), we get
	\begin{align*}
		\left\|\tau(t)\right\|\le \left\|\int_{0}^{b}G(t,s)f(s)\,ds\right\|+\left\|\int_{0}^{b}G(t,s)\mathrm{B}u(s)\mathrm{d}s\right\|.
	\end{align*}
	
	Let us consider
	\begin{align*}
		&M_{\mathrm{B}}=\|\mathrm{B}\|, \quad \beta=\left(\frac{1-\gamma}{\alpha-\gamma} b^{(\alpha-\gamma) /(1-\gamma)}\right)^{1-\gamma},\quad 
		\Lambda_{0}=\frac{\sum _{k=1}^{m} c_{k}}{1-\frac{M}{\Gamma(\alpha)}\sum_{k=1}^{m}\left |c_{k}{t_{k}}^{\alpha-1}\right |},\\
		&\quad \Lambda_{1}=\frac{M}{\Gamma(\alpha)}\beta\left(b^{\alpha-1}M\Lambda_{0}+1\right),\quad \Lambda_{2}=\left(\Lambda_{0}b^{\alpha-1}\frac{M}{\Gamma(\alpha)+1}\right).
	\end{align*}
	
	\begin{align*}
		\left\|\int_{0}^{b}G(t,s)f(s)\,ds\right\|&\le b^{\alpha-1}\frac{M}{\Gamma(\alpha)}\|O\|\int_{0}^{t_k}\sum_{k=1}^{m}\chi_{t_k}(s){\left(t_k-s\right)}^{\alpha-1}\|T_{\alpha}(t_k-s)\|\|f(s)\|\mathrm{d}s\\
		&\quad+\int_{0}^{t}\sum_{k=1}^{m}\left|\chi_{t}(s)\right|{(t-s)}^{\alpha-1}\left \|T_{\alpha}(t-s)\right\|\left\|f(s)\right\|\mathrm{d}s,\\
		&\le \frac{b^{\alpha-1}M^{2}}{\Gamma(\alpha)}\Lambda_{0}\int_{0}^{t_{k}}{(t_{k}-s)}^{\alpha-1} P(s)\psi(r)\mathrm{d}s\\
		&\quad+\frac{M}{\Gamma(\alpha)}\int_{0}^{t}{(t-s)}^{\alpha-1} P(s)\psi(r)\mathrm{d}s,\\
		&\le \frac{M}{\Gamma(\alpha)}\psi(r){\left\|P\right\|}_{L_{\frac{1}{\gamma}}}{\left(\frac{1-\gamma}{\alpha-\gamma }b^{\frac{\alpha-\gamma}{1-\gamma}}\right)}^{1-\gamma}\left(b^{\alpha-1}M\Lambda_{0}+1\right),\\
		&=\frac{M}{\Gamma(\alpha)}\psi(r){\left\|P\right\|}_{L_{\frac{1}{\gamma}}}\beta\left(b^{\alpha-1}M\Lambda_{0}+1\right),\\&=\psi(r){\left\|P\right\|}_{L_{\frac{1}{\gamma}}}\Lambda_{1}.
	\end{align*}
	We have the norm of $G$ as
	\begin{align*}
		\left\|G(t,s)\right\|=\frac{M}{\Gamma(\alpha)}{\left(b-s\right)}\Lambda_{2},
	\end{align*}
	and of $u$ as
	\begin{align*}
		\left\|u(s)\right\|=\frac{M_{\mathrm{B}}M}{a\Gamma(\alpha)}\Lambda_{2}\left(\left\|x_{1}\right\|+\psi(r)+{\left\|P\right\|}_{L_{\frac{1}{\gamma}}}\Lambda_{1}\right),
	\end{align*}
	and
	\begin{align*}
		\left\|\int_{0}^{b}G(t,s)\mathrm{B}u(s)\mathrm{d}s\right\|&\le \frac{{M_{\mathrm{B}}}^{2}M^{2}\Lambda_{1}\Lambda_{2}}{a{\Gamma(\alpha)}^{2}}\left(\left\|x_{1}+\right\|\psi(r){\left\|P\right\|}_{L_{\frac{1}{\gamma}}}\Lambda_{1}\right)\int_{0}^{t}{\left(t-s\right)}^{2\alpha-2}\mathrm{d}s,\\
		&\le\frac{{M_{\mathrm{B}}}^{2}M^{2}{b}^{2\alpha-1}\Lambda_{1}\Lambda_{2}}{a{\Gamma(\alpha)}^{2}\left(2\alpha-1\right)}\left(\left\|x_{1}+\right\|\psi(r){\left\|P\right\|}_{L_{\frac{1}{\gamma}}}\Lambda_{1}\right).
	\end{align*}
	Now we have,
	\begin{align*}
		t^{1-\alpha}\left\|\tau(t)\right\|&=t^{1-\alpha}\left\|\int_{0}^{b}G(t,s)\left[f(s)+\mathrm{B}u(s)\right]\mathrm{d}s\right\|,\\
		&\le t^{1-\alpha}\left\|\int_{0}^{b}G(t,s)f(s)\,ds\right\|+t^{1-\alpha}\left\|\int_{0}^{b}G(t,s)Bu(s)\mathrm{d}s\right\|,\\
		&\le b^{1-\alpha}\frac{M}{\Gamma(\alpha)}\psi(r){\left\|P\right\|}_{L_{\frac{1}{\gamma}}}\beta\left(b^{\alpha-1}M\Lambda_{0}+1\right) \\
		&\quad+\frac{{M_{\mathrm{B}}}^{2}M^{2}{b}^{2\alpha-1}\Lambda_{1}\Lambda_{2}}{a{\Gamma(\alpha)}^{2}\left(2\alpha-1\right)}\left(\left\|x_{1}+\right\|\psi(r){\left\|P\right\|}_{L_{\frac{1}{\gamma}}}\Lambda_{1}\right).
	\end{align*}
	Thus,
	\begin{align*}
		r&\le b^{1-\alpha}\frac{M}{\Gamma(\alpha)}\psi(r){\left\|P\right\|}_{L_{\frac{1}{\gamma}}}\beta\left(b^{\alpha-1}M\Lambda_{0}+1\right)
		+\frac{{M_{\mathrm{B}}}^{2}M^{2}{b}^{\alpha}\Lambda_{1}\Lambda_{2}}{a{\Gamma(\alpha)}^{2}\left(2\alpha-1\right)}\left(\left\|x_{1}\right\|+\psi(r){\left\|P\right\|}_{L_{\frac{1}{\gamma}}}\Lambda_{1}\right),\\ 
		r&\le  \left\{\frac{M}{\Gamma(\alpha)}\beta\left(M\Lambda_{0}+b^{1-\alpha}\right)+\frac{{M_{\mathrm{B}}}^{2}M^{2}{b}^{\alpha}\Lambda_{1}\Lambda_{2}}{a{\Gamma(\alpha)}^{2}\left(2\alpha-1\right)} \right\}\psi(r){\left\|P\right\|}_{L_{\frac{1}{\gamma}}}+\frac{{M_{\mathrm{B}}}^{2}M^{2}{b}^{\alpha}\Lambda_{1}\Lambda_{2}}{a{\Gamma(\alpha)}^{2}\left(2\alpha-1\right)}\left\|x_{1}\right\|.
	\end{align*}
	Dividing both sides by $r$ and taking the low limit as $r\rightarrowtail\infty $, we get
	\begin{align*}
		1\le \lim_{r \rightarrow \infty}inf\frac{\psi(r)}{r}{\left\|P\right\|}_{L_{\frac{1}{\gamma}}}, 
	\end{align*}
	which is a contradiction to (H6).
	
	\textbf{Step$3$:} $\Phi_{a}(x)$ is closed for each $x\in C_{1-\alpha}(J, \mathbb{H})$.
	
	For each given $x\in C_{1-\alpha}(J, \mathbb{H})$, let $\{\tau_{n}\}_{n\ge 0}\subset \Phi_{a}(x) $ such that $\tau_{n}\to \tau$ in $ C_{1-\alpha}(J, \mathbb{H})$. Then there exists $f_{n}\in S_{\partial F,x}$ such that for all $t\in J$
	\begin{align*}
		\tau_{n}(t)&=\int_{0}^{b} G(t,s)[f_{n}(s)+\mathrm{B}u_{n}(s)]\mathrm{d}s,\\
		\text{where }u_{n}(t)&= \mathrm{B}^{*}G^{*}(b,t) \mathrm{R}\left(a, \Gamma_{0}^{b}\right)\{x_{1}-\int_{0}^{b}G(b,\mu)f_{n}(\mu)\mathrm{d}d\mu\}.
	\end{align*}
	From \cite{papageorgiou1985theory}, Propositions 3.1, $S_{\partial F,x}$ is weakly compact in $L^1(J, \mathbb{H})$
	which implies that $f_{n}$ converges weakly to some $f\in S_{\partial F,x}$ in $L^1(J, \mathbb{H})$. Thus, $u_{n}\rightharpoonup u$ and 
	\begin{align*}
		u_{n}(t)&= \mathrm{B}^{*}G^{*}(b,t) \mathrm{R}\left(a, \Gamma_{0}^{b}\right)\{x_{1}-\int_{0}^{b}G(b,\mu)f_{n}(\mu)\mathrm{d}\mu\}.
	\end{align*}
	Then for each $t\in J, \tau_{n}\to \tau(t)$
	\begin{align*}
		\tau(t)=\int_{0}^{b} G(t,s)f(s)\,ds+\int_{0}^{b}G(t,s)\mathrm{B}\mathrm{B}^{*}G^{*}(t,s) \mathrm{R}\left(a, \Gamma_{0}^{b}\right) \{x_{1}-\int_{0}^{b}G(b,\mu)f(\mu)\mathrm{d}\mu\}\mathrm{d}s.
	\end{align*}
	
	Thus we showed the closedness of $\Phi_{a}(x)\forall x\in C_{1-\alpha}(J, \mathbb{H})$.
	
	\textbf{Step$4$:}$\Phi_{a}$ is upper semicontinuous and condensing.
	
	We have
	\begin{align*}
		\Phi_{a}(x)=\{\tau \in C_{1-\alpha}(J, \mathbb{H}): \tau(t)=\int_{0}^{b}G(t,s)[f(s)+Bu(s)]\mathrm{d}s, \quad f \in S_{\partial F, x}, t \in\left(0,b\right].
	\end{align*}
	Now we prove $\Phi_{a}(x)$ is upper semicontinuous and completely continuous. We subdivide the proof into several claims.
	
	\textbf{Claim $1$:} There exists a $r>0$ such that $\Phi_{a}(B_r)\subseteq B_r$.
	
	By utilizing the method employed in step 2, it becomes straightforward to demonstrate the existence of $r>0$ such that $\Phi_{a}(B_r)\subseteq B_r$. 
	
	\textbf{Claim $2$:} $\Phi_{a}(B_r)$ is a family of equicontinuous function
	
	Let $0\le s\le t_1\le t_2\le b$. For each $x\in B_r$, $\phi\in\Phi_{a}(x) $, $\exists f\in S_{\partial F,x}$ such that
	\begin{align*}
		\tau(t)=\int_{0}^{b}G(t,s)[f(s)+\mathrm{B}u(s)]\mathrm{d}s,
	\end{align*}
	then we have
	\begin{align*}
		\tau_2-\tau_1
		&=\int_{0}^{b}\left[G(t_2,s)-G(t_1,s)\right][f(s)+\mathrm{B}u(s)]\mathrm{d}s,\\
		&=\int_{0}^{b}\bigg[\sum_{k=1}^{m}\chi_{t_{k}}{t_{2}}^{\alpha-1}T_{\alpha}(t_{2})O{(t_{k}-s)}^{\alpha-1}T_{\alpha}(t_{k}-s)\\&+\chi_{t_{2}}(s){(t_{2}-s)}^{\alpha-1}T_{\alpha}(t_{2}-s)\bigg][f(s)+\mathrm{B}u(s)]\mathrm{d}s\\
		&-\int_{0}^{b}\bigg[\sum_{k=1}^{m}\chi_{t_{k}}{t_{1}}^{\alpha-1}T_{\alpha}(t_{1})O{(t_{k}-s)}^{\alpha-1}T_{\alpha}(t_{k}-s)\\&+\chi_{t_{1}}(s){(t_{1}-s)}^{\alpha-1}T_{\alpha}(t_{1}-s)\bigg][f(s)+\mathrm{B}u(s)]\mathrm{d}s.
	\end{align*}
	Now,
	\begin{align*}
		&\|\tau_2-\tau_1\|\\&\le\left\|[{t_2}^{\alpha-1}T_{\alpha}(t_2)-{t_1}^{\alpha-1}T_{\alpha}(t_1)]\int_{0}^{b}\sum_{k=1}^{m}\chi_{t_{k}}(s)O{(t_k-s)}^{\alpha-1}T_{\alpha}(t_k-s)\left[f(s)+\mathrm{B}u(s)\right]\mathrm{d}s\right\|\\
		&+\left\|\int_{0}^{b}\chi_{t_{1}}(s)[{(t_2-s)}^{\alpha-1}T_{\alpha}(t_2-s)-{(t_1-s)}^{\alpha-1}T_{\alpha}(t_1-s)][f(s)+\mathrm{B}u(s)]\mathrm{d}s\right\| \\
		&+\left\|\int_{0}^{b}\left[\chi_{t_2}-\chi_{t_1}(s)\right]{(t_2-s)}^{\alpha-1}T_{\alpha}(t_2-s)[f(s)+\mathrm{B}u(s)]\mathrm{d}s\right\|,\\
		&\le\|\left[{t_2}^{\alpha-1}T_{\alpha}(t_2)-{t_1}^{\alpha-1}T_{\alpha}(t_1)\right]\int_{0}^{b}\sum_{k=1}^{m}\chi_{t_{k}}O{(t_{k}-s)}^{\alpha-1}T_{\alpha}(t_{k}-s)[f(s)+\mathrm{B}u(s)]\mathrm{d}s\\
		&+\int_{0}^{t_1}\left\|\left[{(t_{2}-s)}^{\alpha-1}T_{\alpha}(t_{2}-s)-{(t_{1}-s)}^{\alpha-1}T_{\alpha}(t_{1}-s)\right][f(s)+\mathrm{B}u(s)]\mathrm{d}s\right\|\\
		&+\int_{t_1}^{t_2}\left\|{(t_{2}-s)}^{\alpha-1}T_{\alpha}(t_{2}-s)[f(s)+\mathrm{B}u(s)]\mathrm{d}s\right\|,\\	&=\int_{0}^{b}\bigg\|\left[{t_{2}}^{\alpha-1}\{T_{\alpha}(t_{2})-T_{\alpha}(t_{1})\}+\left({t_{2}}^{\alpha-1}-{t_{1}}^{\alpha-1}T_{\alpha}(t_{1})\right)\right]\\&\qquad \hspace{2em}\sum_{k=1}^{m}\chi_{t_{k}}O{(t_k-s)}^{\alpha-1}T_{\alpha}(t_{k}-s)[f(s)+\mathrm{B}u(s)]\mathrm{d}s\bigg\|\\
		&+\int_{0}^{t_{1}}\bigg\|\left[{(t_{2}-s)}^{\alpha-1}\left\{T_{\alpha}(t_{2}-s)-T_{\alpha}(t_{1}-s)\right\}+\left\{{(t_{2}-s)}^{\alpha-1}-{(t_{1}-s)}^{\alpha-1}\right\}T_{\alpha}(t_{1}-s)\right]\\&
		\qquad \hspace{2em} [f(s)+\mathrm{B}u(s)]\mathrm{d}s\bigg\|+\int_{t_1}^{t_2}\left\|{(t_{2}-s)}^{\alpha-1}T_{\alpha}(t_{2}-s)[f(s)+\mathrm{B}u(s)]\mathrm{d}s\right\|,\\
		&\le \Lambda_{0}\max_{t_1,t_2\in [0,b]}\{T_{\alpha}(t_2)-T_{\alpha}(t_1)\}b^{\alpha-1}\frac{M}{\Gamma_{\alpha}}\bigg[\int_{0}^{t_k}{(t_{k}-s)}^{\alpha-1}\|\mathrm{B}B\|\|u(s)\|\,ds\\&\qquad+\int_{0}^{t_k}{(t_{k}-s)}^{\alpha-1}\|f(s)\|\mathrm{d}s\bigg]\\
		&+\max_{s\in[0,t_{1}]}\left\|T_{\alpha}(t_{2}-s)-T_{\alpha}(t_{1}-s)\right\|\int_{0}^{t_1}{(t_{2}-s)}^{\alpha-1}\left\|[f(s)+\mathrm{B}u(s)]\right\|\mathrm{d}s\\
		&+\frac{M}{\Gamma_{\alpha}}\int_{0}^{t_1}\{{(t_{2}-s)}^{\alpha-1}-{(t_{1}-s)}^{\alpha-1}\}\|[f(s)+\mathrm{B}u(s)]\|\,ds\\&+\frac{M}{\Gamma_{\alpha}}\int_{t_1}^{t_2}{(t_{2}-s)}^{\alpha-1}\|[f(s)+\mathrm{B}u(s)]\|\mathrm{d}s,\\
		&\le \Lambda_{0}\max_{t_1,t_2\in [0,b]}\{T_{\alpha}(t_2)-T_{\alpha}(t_1)\}b^{\alpha-1}\frac{MM_\mathrm{B}}{\Gamma_{\alpha}}\int_{0}^{t_k}{(t_{k}-s)}^{\alpha-1}\|u(s)\|\mathrm{d}s\\
		&+\Lambda_{0}\max_{t_1,t_2\in[0,b]}\{T_{\alpha}(t_2)-T_{\alpha}(t_1)\}b^{\alpha-1}\frac{MM_\mathrm{B}}{\Gamma_{\alpha}}\int_{0}^{t_k}{(t_{k}-s)}^{\alpha-1}\|f(s)\|\mathrm{d}s\\
		&+M_{\mathrm{B}}\max_{s\in[0,t_{1}]}\left\|T_{\alpha}(t_{2}-s)-T_{\alpha}(t_{1}-s)\right\|\int_{0}^{t_1}{(t_{2}-s)}^{\alpha-1}\|u(s)\|\mathrm{d}s\\
		&+\max_{s\in[0,t_{1}]}\left\|T_{\alpha}(t_{2}-s)-T_{\alpha}(t_{1}-s)\right\|\int_{0}^{t_1}{(t_{2}-s)}^{\alpha-1}\|f(s)\|\mathrm{d}s\\
		&+\frac{MM_{\mathrm{B}}}{\Gamma_{\alpha}}\int_{0}^{t_1}\{{(t_{2}-s)}^{\alpha-1}-{(t_{1}-s)}^{\alpha-1}\}\|u(s)\|\mathrm{d}s\\
		&+\frac{M}{\Gamma_{\alpha}}\int_{0}^{t_1}\{{(t_{2}-s)}^{\alpha-1}-{(t_{1}-s)}^{\alpha-1}\}\|f(s)\|\mathrm{d}s\\
		&+\frac{MM_{\mathrm{B}}}{\Gamma_{\alpha}}\int_{t_1}^{t_2}{(t_{2}-s)}^{\alpha-1}\|u(s)\|\,ds+\frac{M}{\Gamma_{\alpha}}\int_{t_1}^{t_2}{(t_{2}-s)}^{\alpha-1}\|f(s)\|\mathrm{d}s,\\
	\end{align*}
	\begin{align}\label{integraleq}
		\le I_{1}+I_{2}+I_{3}+I_{4}+I_{5}+I_{6}+I_{7}+I_{8},
	\end{align}
	where
	\begin{align*}
		&I_1:=\Lambda_{0}\max_{t_1,t_2\in [0,b]}\{T_{\alpha}(t_2)-T_{\alpha}(t_1)\}b^{\alpha-1}\frac{MM_{\mathrm{B}}}{\Gamma_{\alpha}}\int_{0}^{t_k}{(t_{k}-s)}^{\alpha-1}\|u(s)\|\mathrm{d}s,\\
		&I_2:=\Lambda_{0}\max_{t_1,t_2\in[0,b]}\{T_{\alpha}(t_2)-T_{\alpha}(t_1)\}b^{\alpha-1}\frac{M}{\Gamma_{\alpha}}\|P\|_{\frac{1}{\gamma}}\psi(r)\frac{1-\gamma}{\alpha-\gamma}{t_k}^{\frac{\alpha-\gamma}{1-\gamma}},\\
		&I_3:=M_{\mathrm{B}}\max_{s\in[0,t_{1}]}\left\|T_{\alpha}(t_{2}-s)-T_{\alpha}(t_{1}-s)\right\|\int_{0}^{t_1}{(t_{2}-s)}^{\alpha-1}\|u(s)\|\mathrm{d}s,\\
		&I_4:=M_{\mathrm{B}}\max_{s\in[0,t_{1}]}\left\|T_{\alpha}(t_{2}-s)-T_{\alpha}(t_{1}-s)\right\|\|P\|_{\frac{1}{\gamma}}\psi(r)\left[\left(\frac{1-\gamma}{\alpha-\gamma}\right)\left\{{t_2}^\{\frac{\alpha-\gamma}{1-\gamma}-{t_2-t_1}^{\frac{\alpha-\gamma}{1-\gamma}}\right\}\right],\\
		&I_5:=\frac{M}{\Gamma_{\alpha}}\|P\|_{\frac{1}{\gamma}}\psi(r)\left(\frac{1-\gamma}{\alpha-\gamma}\right)\left[-{(t_2-t_1)}^{\frac{\alpha-\gamma}{1-\gamma}}+{t_2}^{\frac{\alpha-\gamma}{1-\gamma}}-{t_1}^{\frac{\alpha-\gamma}{1-\gamma}}\right],\\
		&I_6:=\frac{MM_{\mathrm{B}}}{\Gamma_{\alpha}}\int_{t_1}^{t_2}{(t_{2}-s)}^{\alpha-1}\|u(s)\|\mathrm{d}s,\\
		&I_7:=\frac{M}{\Gamma_{\alpha}}\|P\|_{\frac{1}{\gamma}}\psi(r)\left(\frac{1-\gamma}{\alpha-\gamma}\right){(t_2-t_1)}^{\frac{\alpha-\gamma}{1-\gamma}},\\
		&I_8:=\frac{MM_{\mathrm{B}}}{\Gamma_{\alpha}}\int_{t_1}^{t_2}{(t_{2}-s)}^{\alpha-1}\|u(s)\|\mathrm{d}s.
	\end{align*}
	From lemma 2.1 $T_{\alpha}(t)$ is continuous in the uniform operator topology for $t>0$. From this property of $T_{\alpha}$ we directly obtain $I_{1}$, $I_{2}$, $I_{3}$, $I_{4}$ tends to $0$ independently of $x\in B_{r}$ as $t_{2}\rightarrow t_{1}$. $I_{5}$ also tends to $0$ independently of $x\in B_{r}$ as $t_{2}\rightarrow t_{1}$. Using the absolute continuity of Lebesgue integral, we have $I_{6}$, $I_{7}$, $I_{8}$ tending to $0$ independently of $x\in B_{r}$ as $t_{2}\rightarrow t_{1}$.
	
	Therefore, $\Phi_{a}(B_r)\subset C_{1-\alpha}(J, \mathbb{H})$ is a family of equicontinuous function.
	
	\textbf{Claim $3$:} The set $\Pi(t)=\{\tau(t):\tau\in\Phi_{a}(B_r)\}\subset\mathbb{H}$ is relatively compact for each $t\in J$.
	
	Let $0<t\le b$ be fixed. For $x\in B_{r}$ and $\tau\in\Phi_{a}(x)$, $f\in S_{\partial F,x}$ such that for each $t\in J$,
	\begin{align*}
		\tau(t)&=\sum_{k=1}^{m}\int_{0}^{t_k}c_{k}t^{\alpha-1}T_{\alpha}(t)O{(t_{k}-s)}^{\alpha-1}T_{\alpha}(t_{k}-s)[f(s)+\mathrm{B}u(s)]\mathrm{d}s\\&\quad+\int_{0}^{t}{(t-s)}^{\alpha-1}T_{\alpha}(t-s)[f(s)+\mathrm{B}u(s)]\mathrm{d}s,
	\end{align*}
	where
	\begin{align*}
		u(t)=\mathrm{B}^{*}&\bigg(\sum_{k=1}^{m}\chi_{t_k}(s)t^{\alpha-1}T_{\alpha}^{*}(t)o^{*}{(t_{k}-s)}^{\alpha-1}T_{\alpha}^{*}(t_{k}-s)+\chi_{t}(s){(t-s)}^{\alpha-1}T_{\alpha}^{*}(t-s)\bigg)R(a,\Gamma_{0}^{b})\\&\bigg(x_{1}-\sum_{k=1}^{m}\int_{0}^{t_k}c_{k}t^{\alpha-1}T_{\alpha}(b)o{(t_{k}-s)}^{\alpha-1}T_{\alpha}^{*}(t_{k}-s)f(s)\mathrm{d}s\\&\quad-\int_{0}^{b}{(b-s)}^{\alpha-1}T_{\alpha}(b-s)f(s)\mathrm{d}s\bigg).
	\end{align*}
	For all $\epsilon\in(0,t)$ and for all $\delta>0$, define
	\begin{align*}
		&\tau^{\epsilon,\delta}(t)\\&=\alpha t^{\alpha-1}\mathcal{T}_{\alpha}(t)\mathcal{T}(\epsilon^{\alpha}\delta)\sum_{k=1}^{m}c_{k}o\int_{0}^{t_{k}-\epsilon}\int_{\delta}^{\infty}\theta{(t_{k}-s)}^{\alpha-1}\xi_{\alpha}(\theta)\mathcal{T}({(t_{k}-s)}^{\alpha}\theta-\epsilon^{\alpha}\delta)f(s)\mathrm{d}\theta\mathrm{d}s\\
		&\quad+\alpha t^{\alpha-1}\mathcal{T}(t)\mathcal{T}(\epsilon^{\alpha}\delta)\sum_{k=1}^{m}c_{k}o\int_{0}^{t_{k}-\epsilon}\int_{\delta}^{\infty}\theta{(t_{k}-s)}^{\alpha-1}\xi_{\alpha}(\theta)\mathcal{T}({(t_{k}-s)}^{\alpha}\theta-\epsilon^{\alpha}\delta)\mathrm{B}\mathrm{B}^{*}u(t)\mathrm{d}\theta\mathrm{d}s\\
		&\quad+\alpha\mathcal{T}(\epsilon^{\alpha}\delta)\int_{0}^{t-\epsilon}\int_{\delta}^{\infty}\theta{(t-s)}^{\alpha-1}\xi_{\alpha}(\theta)\mathcal{T}({(t-s)}^{\alpha}\theta-\epsilon^{\alpha}\delta)f(s)\mathrm{d}\theta\mathrm{d}s\\
		&\quad+\alpha\mathcal{T}(\epsilon^{\alpha}\delta)\int_{0}^{t-\epsilon}\int_{\delta}^{\infty}\theta{(t-s)}^{\alpha-1}\xi_{\alpha}(\theta)\mathcal{T}({(t-s)}^{\alpha}\theta-\epsilon^{\alpha}\delta)\mathrm{B}\mathrm{B}^{*}u(t)\mathrm{d}\theta\mathrm{d}s.
	\end{align*}
	By the compactness of $\mathcal{T}(\epsilon^{\alpha}\delta)(\epsilon^{\alpha}\delta>0)$, we obtain the set $\Pi^{\epsilon,\delta}(t)=\{\phi^{\epsilon,\delta}(t);\tau\in \Phi_{a}(B_{r})\}$ which is relatively compact in $\mathbb{H}\forall\epsilon\in (0,t)$ and $\delta>0$, moreover we have
	\begin{align*}
		&\|\tau(t)-\tau^{\epsilon,\delta}\|\\&\le\bigg\|\alpha \mathcal{T}_{\alpha}(t)t^{\alpha-1}\sum_{k=1}^{m}\int_{0}^{t_k}\int_{0}^{\infty}c_{k}o{(t_k-s)}^{\alpha-1}\mathcal{T}({(t_k-s)}^{\alpha}\theta)\xi_{\alpha}(\theta)\theta [f(s)+\mathrm{B}u(s)]\mathrm{d}\theta\mathrm{d}s\\
		&\quad+\alpha \int_{0}^{t}\int_{0}^{\infty}\theta {(t-s)}^{\alpha-1}\xi_{\alpha}(\theta)\mathcal{T}({(t-s)}^{\alpha}\theta))[f(s)+\mathrm{B}u(s)]\mathrm{d}\theta\mathrm{d}s\\
		&\quad+\alpha \mathcal{T}_{\alpha}(t)t^{\alpha-1}\mathcal{T}(\epsilon^{\alpha}\delta)\sum_{k=1}^{m}\int_{0}^{t_k-\epsilon}\int_{\delta}^{\infty}c_{k}o{(t_k-s)}^{\alpha-1}\mathcal{T}({(t_k-s)}^{\alpha}\theta-\epsilon^{\alpha}\delta)\xi_{\alpha}(\theta)\theta\\&\qquad [f(s)+\mathrm{B}u(s)]\mathrm{d}\theta\mathrm{d}s\\
		&\quad+\alpha\mathcal{T}(\epsilon^{\alpha}\delta) \int_{0}^{t-\epsilon}\int_{\delta}^{\infty}\theta {(t-s)}^{\alpha-1}\xi_{\alpha}(\theta)\mathcal{T}({(t-s)}^{\alpha}\theta-\epsilon^{\alpha}\delta)[f(s)+\mathrm{B}u(s)]\mathrm{d}\theta\mathrm{d}s\bigg\|,\\
		&=\bigg\|\alpha \mathcal{T}_{\alpha}(t)t^{\alpha-1}\bigg[\sum_{k=1}^{m}\int_{0}^{t_k}\int_{0}^{\infty}c_{k}o{(t_k-s)}^{\alpha-1}\mathcal{T}({(t_k-s)}^{\alpha}\theta)\xi_{\alpha}(\theta)\theta [f(s)+\mathrm{B}u(s)]\mathrm{d}\theta\mathrm{d}s\\
		&\quad-\sum_{k=1}^{m}\int_{0}^{t_k-\epsilon}\int_{\delta}^{\infty}c_{k}o{(t_k-s)}^{\alpha-1}\mathcal{T}({(t_k-s)}^{\alpha}\theta)\xi_{\alpha}(\theta)\theta [f(s)+\mathrm{B}u(s)]\mathrm{d}\theta\mathrm{d}s\bigg]\\
		&\quad+\alpha\bigg[\int_{0}^{t}\int_{0}^{\infty}\theta {(t-s)}^{\alpha-1}\xi_{\alpha}(\theta)\mathcal{T}({(t-s)}^{\alpha}\theta)[f(s)+\mathrm{B}u(s)]\mathrm{d}\theta\mathrm{d}s\\
		&\quad-\int_{0}^{t-\epsilon}\int_{\delta}^{\infty}\theta {(t-s)}^{\alpha-1}\xi_{\alpha}(\theta)\mathcal{T}({(t-s)}^{\alpha}\theta)[f(s)+\mathrm{B}u(s)]\mathrm{d}\theta\mathrm{d}s\bigg]\bigg\|,\\
		&\le \left\|\alpha \mathcal{T}_{\alpha}t^{\alpha-1}\left[\sum_{k=1}^{m}\int_{t_k-\epsilon}^{t_k}\int_{0}^{\delta}c_{k}o\theta{(t_k-s)}^{\alpha-1}\mathcal{T}({(t_k-s)}^{\alpha-1}\theta)\xi_{\alpha}(\theta)\{f(s)+\mathrm{B}u(s)\}\mathrm{d}\theta\mathrm{d}s \right]\right \|\\ &\quad+\left\|\alpha\left[\int_{t-\epsilon}^{t}\int_{0}^{\delta}\theta {(t-s)}^{\alpha-1}\xi_{\alpha}(\theta)\mathcal{T}({(t-s)}^{\alpha}\theta)\{\mathrm{B}u(s)+f(s)\}\mathrm{d}\theta\mathrm{d}s\right]\right\|,\\
		&\le \alpha \frac{M^2}{\Gamma_{\alpha}}b^{\alpha-1}\Lambda_{0}M_{\mathrm{B}}\int_{t_k-\epsilon}^{t_k}{(t_k-s)}^{\alpha-1}\|u(s)\|\mathrm{d}s\int_{0}^{\delta}\theta \xi_{\alpha}(\theta)\mathrm{d}\theta\\
		&\quad+\alpha\frac{M^2}{\Gamma_{\alpha}}b^{\alpha-1}\Lambda_{0}\|P\|_{\frac{1}{\gamma}}\psi(r)\int_{t_k-\epsilon}^{t_k}{(t_k-s)}^{\alpha-1}\mathrm{d}s\int_{0}^{\delta}\xi_{\alpha}\mathrm{d}\theta\\
		&\quad+\alpha M\|P\|_{\frac{1}{\gamma}}\psi(r)\int_{t-\epsilon}^{t}{(t-s)}^{\alpha -1}\mathrm{d}s\int_{0}^{\delta}\theta \xi_{\alpha}(\theta)\mathrm{d}\theta\\
		&\quad+\alpha M M_{\mathrm{B}}\int_{t-\epsilon}^{t}{(t-s)}^{\alpha -1}\|u(s)\|\mathrm{d}s\int_{0}^{\delta}\theta \xi_{\alpha}(\theta)\mathrm{d}\theta. 
	\end{align*}
	In the above inequality, as $\epsilon$ approaches zero, the right-hand side of the inequality approaches zero as well. This implies that there exist relatively compact sets that are arbitrarily close to the set  $\Pi(t)$ for $t>0$. Consequently, the set $\Pi(t)$, $t>0$ is also relatively compact in $\mathbb{H}$. Combining Claims $1–3$ with the Arzola-Ascoli theorem, we can deduce that $\Phi_{a}$ is a completely continuous function.
	
	\textbf{Claim $4$:} $\Phi_{a}$ has a closed graph.
	
	Let $x_{n}\rightarrow x^{*}(n\rightarrow\infty)$, $\tau_{n}\in \Phi_{a}(x_n)$, $\tau_{n}\rightarrow \tau^{*}(n\rightarrow \infty)$. Our aim is to prove $\tau^{*}\in \Phi_{a}(x^{*})$. Since $\tau_{n}\in \Phi_{a}(x_n)$, $\exists f_{n}\in S_{\partial F,x_{n}}$ such that for each $t\in J$ we have
	\begin{align*}
		\tau_{n}(t)&=\int_{0}^{b} G(t,s)f_{n}(s)\mathrm{d}s+\int_{0}^{b}G(t,s)\mathrm{B}\mathrm{B}^{*}G^{*}(t,s)\mathrm{R}\left(a, \Gamma_{0}^{b}\right)\left\{x_{1}-\int_{0}^{b}G(b,\mu)f_{n}(\mu)\mathrm{d}\mu\right\}\mathrm{d}s,\\
		&=\sum_{k=1}^{m}c_k\int_{0}^{t_k}t^{\alpha-1}\mathcal{T_\alpha}(t)\mathcal{O}{(t_k-s)}^{\alpha-1}\mathcal{T_\alpha}(t_k-s)f_n(s)\mathrm{d}s+\int_{0}^{t}{(t-s)}^{\alpha-1}\mathcal{T_\alpha}(t-s)f_n(s)\mathrm{d}s\\
		&\quad+\int_{0}^{b}G(t,s)\mathrm{B}\mathrm{B}^{*}G^{*}(t,s)\mathrm{R}\left(a, \Gamma_{0}^{b}\right)\left\{x_1-\int_{0}^{b}G(b,s)f_n(\mu )\mathrm{d}\mu\right\}\mathrm{d}s.
	\end{align*}
	we must prove that $\exists$ $f^{*}(s)\in S_{\partial F,x^{*}}$, such that $\forall t\in J$,
	\begin{align*}
		\tau^{*}(t)=\int_{0}^{b}G(t,s)f^{*}(s)\mathrm{d}s+\int_{0}^{b}G(t,s)\mathrm{B}\mathrm{B}^{*}G^{*}(t,s)\mathrm{R}\left(a, \Gamma_{0}^{b}\right)\left\{x_1-\int_{0}^{b}G(b,s)f^{*}(\mu)\mathrm{d}\mu\right\}\mathrm{d}s.
	\end{align*}
	Since $\tau_{n}\rightarrow\tau^{*}(n\rightarrow\infty)$, we can obtain
	\begin{align*}
		&\|\int_{0}^{b}G(t,s)f_n(s)\mathrm{d}s+\int_{0}^{b}G(t,s)\mathrm{B}\mathrm{B}^{*}G^{*}(t,s)\mathrm{R}\left(a, \Gamma_{0}^{b}\right)\left\{x_1-\int_{0}^{b}G(b,s)f_n(\mu)\mathrm{d}\mu\right\}\mathrm{d}s\\
		&-\int_{0}^{b}G(t,s)f^{*}(s)\mathrm{d}s+\int_{0}^{b}G(t,s)\mathrm{B}\mathrm{B}^{*}G^{*}(t,s)\mathrm{R}\left(a, \Gamma_{0}^{b}\right)\left\{x_1-\int_{0}^{b}G(b,s)f^{*}(\mu)\mathrm{d}\mu\right\}\mathrm{d}s\| \rightarrow 0 \text{ as } n\rightarrow \infty.
	\end{align*}
	Consider the linear continuous operator $\Gamma: L^{\frac{1}{\gamma}}(J,\mathbb{H})\rightarrow C_{1-\alpha}(J, \mathbb{H})$
	\begin{align*}
		(\Gamma f)(t)=\int_{0}^{b}G(t,s)f(s)\mathrm{d}s-\int_{0}^{b}G(t,s)\mathrm{B}\mathrm{B}^{*}G^{*}(t,s)\mathrm{R}\left(a, \Gamma_{0}^{b}\right)\left(\int_{0}^{b}G(t,\mu)f(\mu)\mathrm{d}\mu\right)\mathrm{d}s.
	\end{align*}                     
	Clearly it follows from lemma\ref{lem3.2} that $\Gamma o S_{\partial F}$ is a closed graph operator. Moreover, we have
	\begin{align*}
		\tau_{n}(t)-\int_{0}^{b}G(t,s)\mathrm{B}\mathrm{B}^{*}G^{*}(t,s)\mathrm{R}\left(a, \Gamma_{0}^{b}\right)x_{1}\mathrm{d}s\in \Gamma(S_{\partial f,x_{n}}).
	\end{align*}
	Since $x_{n}\longrightarrow x_{*}$, it follows from lemma\ref{lem3.2} that
	\begin{align*}
		\tau_{*}(t)-\int_{0}^{b}G(t,s)\mathrm{B}\mathrm{B}^{*}G^{*}(t,s)\mathrm{R}\left(a, \Gamma_{0}^{b}\right)x_{1}\mathrm{d}s\in \Gamma(S_{\partial f,x_{*}}). 
	\end{align*}
	Therefore $\Phi_{a}$ has a closed graph from lemma\ref{lem3.2}. Since $\Phi_{a}$ is completely continuous multivalued map with compact value, we have that $\Phi_{a}$ is upper semi continuous.
	
	Thus $\Phi_{a}$ is upper semicontinuous and condensing. Therefore by theorem\ref{thrm3.1}, we conclude that $\Phi_{a}$ has a fixed point $x(.)$ on $B_{r_{0}}$. Thus, the fractional control system\eqref{P} has a mild solution.
	\section{\textbf{Approximate Controllability Results}}
	\vspace{0.5cm}
	In this section we obtain sufficient conditions of approximate controllability of the system \eqref{P2}. Motivation is from the case of linear system. Here we additonally assume
	\begin{Ass}
		\begin{itemize}
			\item[(H5$'$)] There exists a positive constant $L$ such that $\|\partial F(t, x(t))\| \leq L$ for all $(t, x) \in J \times \mathbb{H}$.
		\end{itemize}
	\end{Ass}
	
	\begin{theorem}\label{thrm4.1}
		Assume that assumptions (H1)-(H5) and (H5$'$) are satisfied and the linear system \eqref{LP} is approximately controllable on $J$. Then system \eqref{P2} is approximately controllable on $J$.
	\end{theorem}
	Proof: Let $x^{a}$ be a fixed point of $\Phi_{a}$ in $B_{r_0}$, this means that $\exists f^{a}\in S_{\partial F,x^{a}}$ such that $\forall t\in J$,
	\begin{align*}
		x^{a}(t)=\int_{0}^{b}G(t,s)[f^{a}(s)+\mathrm{B}\mathrm{B}^{*}G^{*}(t,s) \mathrm{R}\left(a, \Gamma_{0}^{b}\right)\{x_{1}-\int_{0}^{b}G(b,\mu)f^{a}(\mu) \mathrm{d}\mu\}]\mathrm{d}s.
	\end{align*} 
	Now we define a function
	\begin{align*}
		P(f^{a})=x_{1}-\int_{0}^{b}G(b,s)f^{a}(s)\mathrm{d}s,\quad \text{for some } f^{a}\in S_{\partial F,x^{a}}.
	\end{align*}	
	Note that $ I-\Gamma_{0}^{b}\mathrm{R}\left(a, \Gamma_{0}^{b}\right)=a\mathrm{R}\left(a, \Gamma_{0}^{b}\right)$, we get
	\begin{align*}
		x^{a}(b)=X_{1}-a\mathrm{R}\left(a, \Gamma_{0}^{b}\right)P(f^{a}).
	\end{align*}
	By assumption (H5$'$),
	\begin{align*}
		\int_{0}^{b}\left\|f^{a}(s)\right\|^{2}\mathrm{d}s\le L^{2}b.
	\end{align*}
	This implies that the sequence $\{f^{a}\}$, that converges weakly to say, $f$ in $L^{\frac{1}{\gamma}}(J,\mathbb{H})$. Let us denote
	\begin{align*}
		h=x_{1}-\int_{0}^{b}G(b,s)f(s)\mathrm{d}s,
	\end{align*}	
	we see that 
	\begin{align}\label{eq3.13}
		\left\|P(f^{a})-h\right\|&=\left\|x_{1}-\int_{0}^{b}G(b,s)f^{a}(s)\mathrm{d}s-x_{1}+\int_{0}^{b}G(b,s)f(s)\mathrm{d}s\right\|,\nonumber\\
		&\le \sup_{t \in J}\left\|\int_{0}^{b}G(b,s)\left[f^{a}(s)-f(s)\right]\mathrm{d}s\right\|. 
	\end{align}
	By (H6$'$) and Ascoli- Arzela theorem we can show that the linear operator $g\longrightarrow\int_{0}^{.}G(.,s)g(s)\mathrm{d}s : L^{\frac{1}{\gamma}(J,\mathbb{H})}\longrightarrow  C_{1-\alpha}(J, \mathbb{H})$ is compact, consequently the right hand side of \eqref{eq3.13} tends to zero as $a\longrightarrow0^{+}$. Now
	\begin{align*}
		\left\|x^{a}(b)-x_{1}\right\|&=\left\|a\mathrm{R}\left(a, \Gamma_{0}^{b}\right)P(f^{a})\right\|,\\
		&\le \left\|a\mathrm{R}\left(a, \Gamma_{0}^{b}\right)(h)\right\|+\left\|a\mathrm{R}\left(a, \Gamma_{0}^{b}\right)\left(P(f^{a})-h\right)\right\|,\\
		&\le \left\|a\mathrm{R}\left(a, \Gamma_{0}^{b}\right)(h)\right\|+\left\|\left(P(f^{a})-h\right)\right\|\longrightarrow0, 
	\end{align*}
	as $a\longrightarrow 0^{+}$.
	Thid proves the approximate controllability of system \eqref{P}.
	\section{\textbf{Application}}\setcounter{equation}{0}
	
	In this section, we provide a examples to validate the results obtained in the previous sections.
	
	\begin{Ex}\label{ex1} Let us  consider the following heat conduction system:
		
		\begin{equation}\label{5.1}
			\left\{
			\begin{aligned}
				^R\mathrm{D_{0,t}^{\frac{3}{4}}} x(t,y)-\frac{\partial^{2}}{\partial y^{2}}x(t,y)&=b(y)u(t)+Q(t,y),\quad 0<y<\pi,\quad t\in J=[0,b],\\
				x(t,0)&=x(t,\pi)=0, \quad t\in J,\\
				{I_{0,t}}^{1-\alpha}_{0^{+}}x(t)|_{t=0}&=\sum _{k=1}^{m}c_{k}x(t_k,y),\quad y\in[0,\pi].
			\end{aligned}
			\right.
		\end{equation}
		where $x(t,y)$ represents the temperature at point $y\in[0,\pi]$ and time $t\in J$. $^R\mathrm{D_{0,t}^{\frac{3}{4}}}$ is the R-L fractional derivative of order $\frac{3}{4}$, $J=[0,b]$. It is supposed that $ Q=\bar{Q}+\bar{\bar{Q}}$, where $\bar{\bar{Q}}$ is a continuous function and $\bar{Q}$ is a known function of the temperature of the form $$-\stackrel{-}{Q}\in\partial F(t,x(t,y))\quad (t,y)\in J\times(0,\pi),$$ with a measurable function $F$ provided $F(t,.)$ is locally Lipschitz on $\mathbb{R}$, so its generalized gradient $\partial F$is well defined. For $k=1,2,...,m$ all $c_{k}\in \mathbb{R}$ and satisfy \ref{ass1}.
		
		Let us take $\mathbb{H}=\mathrm{L}^{2}([0,\pi];\mathbb{R})$, and the family of operators $\mathrm{A}$ as $$\mathrm{A}x=\frac{\partial^{2}}{\partial y^{2}}x(t,y),$$ with the domain $\mathrm{D}(\mathrm{A})=\{x\in \mathbb{H}; x, x^{\prime} \text{are absolutely continuous}, x^{\prime\prime}\in\mathbb{H}, x(0)=x(\pi)=0\}$.
		
		Then $$ \mathrm{A}x=-\sum_{n=1}^{\infty}n^{2}\langle x,e_{n}\rangle e_{n}, \quad x\in \mathrm{D}(\mathrm{A}),$$ where $$e_{n}(y)=\sqrt{\frac{2}{\pi}} \sin ny, \quad y\in [0,\pi ], \quad n=1,2,.....,$$ is orthogonal set of eigenvectors of $\mathrm{A}$. It is well known that the operator $\mathrm{A}$ generates a strongly continuous semigroup  $\mathcal{T}(t)(t\ge0)$ on $\mathbb{H}$, which are compact and is given by 
		\begin{align*}
			\mathcal{T}(t)x=\sum_{n=1}^{\infty}e^{-n^{2}t}\langle x,e_{n}\rangle e_{n}, \qquad x\in \mathbb{H}
		\end{align*}
		and \begin{align*}
			\mathcal{T}_{\frac{3}{4}}(t)&=\frac{3}{4}\int_{0}^{\infty}\theta\xi_{\frac{3}{4}}(\theta)\mathcal{T}(t^{\frac{3}{4}}\theta)\mathrm{d}\theta,\\
			\mathcal{T}_{\frac{3}{4}}(t)&=\frac{3}{4}\sum_{n=1}^{\infty}\int_{0}^{\infty}\theta\xi_{\frac{3}{4}}(\theta) \exp(-n^{2}t^{\frac{3}{4}}\theta )\mathrm{d}\theta \langle x,e_{n}\rangle.
		\end{align*}
		$\mathcal{T}_{\frac{3}{4}}$ and $F(t,y)$ satisfy (H2)-(H5).
		
		Let $\mathrm{B}\in L(\mathbb{R},\mathbb{H})$ be defined as,
		$$
		(\mathrm{B}u)(y)=b(y)u,\quad  \; \;\mathrm{B}^{*}v=\sum_{n=1}^{\infty}\langle b, e_{n}\rangle\langle v,e_{n}\rangle,
		$$
		where $y\in [0,\pi]$, $u\in \mathbb{R}$ and $b(y)\in L_{2}[0,\pi]$.

		In order to show that associated linear system is approximate controllable on $[0,b]$, we need to show that $(b-s)^{\alpha-1}\mathrm{B}^{*}\mathcal{T}_{\alpha}(b-s)x=0 \implies x=0$. We observe that for $\frac{1}{2}<\mu \le 1$
		\begin{align}
			(b-s)^{\mu-1}\mathrm{B}^{*}\mathcal{T}_{\alpha}(b-s)x&=(b-s)^{\mu-1}\sum_{n=1}^{\infty}\langle b, e_{n}\rangle \frac{3}{4}\int_{0}^{\infty}\theta \xi_{\frac{3}{4}}(\theta) \exp(-n^{2}t^{\frac{3}{4}}\theta)\mathrm{d}\theta\langle x,e_{n}\rangle,\\
			&=(b-s)^{\mu-1}\frac{3}{4}\sum_{n=1}^{\infty}\int_{0}^{\infty}\theta \xi_{\frac{3}{4}}(\theta) \exp(-n^{2}t^{\frac{3}{4}}\theta)\mathrm{d}\theta\langle b,e_{n}\rangle \langle x,e_{n}\rangle = 0.
		\end{align}
		This gives $\langle x,e_{n}\rangle = 0 \implies x=0 $ provided that $\langle b,e_{n}\rangle=\int_{0}^{\pi }b(\theta)e_{n}\theta\mathrm{d}\theta\neq 0 $ for $n=1,2,...$.
		
		Therefore, the associated linear system is approximate controllable provided that $\int_{0}^{\pi} b(\theta)e_{n}(\theta)\mathrm{d}\theta \neq 0$ for $n=1,2,3,....$. Because of the compactness of the semigroup $\mathcal{T}$ generated by $\mathrm{A}$, the associated linear system is not exactly controllable but it is approximate controllable. Hence from theorem\ref{Thrm3.2} there exists a mild solution of problem\eqref{5.1} and by theorem\ref{thrm4.1} the given system\eqref{5.1} is approximate controllable.
	\end{Ex}
	
	\section{\textbf{Conclusion}}
	This paper explores the existence of mild solutions and approximate controllability for Riemann-Liouville fractional differential Hemivariational inequalities within a separable Hilbert space. Employing nonsmooth analysis and multivalued theory, we utilize fixed-point techniques and ideas from semigroup theory to derive our results. Additionally, we provide an illustrative example to demonstrate the efficacy of our findings. Our future aims include delving into the existence and controllability of Hemivariational Inequality problems within separable reflexive Banach spaces, while also addressing the impulse effect within this framework.
	
	\section{\textbf{Acknowledgment}}
	The corresponding authors are thankful to the funding agency SERB New Delhi for their financial support for project No: MTR/2023/000245.
	
	\section{\textbf{Conflict of Interest}}
	The authors declare that they have no conflict of interest.


\begin{thebibliography}{20}
		
		%
		\bibitem{arora2020approximate} 
		{\sc S. Arora, S. Singh, J. Dabas and  M. T. Mohan},
		{\it Approximate controllability of semilinear impulsive functional differential systems with non-local conditions},
		IMA Journal of Mathematical Control and Information, {\bf 37}, (4) (2020), 1070--1088.
		
		\bibitem{byszewski1999existence}
		{\sc  L. Byszewski},
		{\it Existence and uniqueness of a classical solution to a functional-differential abstract nonlocal Cauchy problem},
		Journal of Applied Mathematics and Stochastic Analysis, {\bf 12}, (1) (1999), 91--97.
		
		\bibitem{chen2014existence}
		{\sc  P. Chen and Y. Li},
		{\it Existence of mild solutions for fractional evolution equations with mixed monotone nonlocal conditions},
		Zeitschrift f{\"u}r angewandte Mathematik und Physik, {\bf 65}, (2014), 711--728 .
		
		\bibitem{chen2020existence} 
		{\sc P. Chen, X. Zhang and Y. Li},
		{\it Existence and approximate controllability of fractional evolution equations with nonlocal conditions via resolvent operators},
		Fractional Calculus and Applied Analysis, {\bf 23}, (1) (2020), 268--291.
		
		\bibitem{clarke1983optimization}
		{\sc F. H. Clarke},
		{\it 	Optimization and Nonsmooth Analysis},
		Wiley, Canadian Mathematical Society series of monographs and advanced texts, (1983).
		
		\bibitem{deimling2011multivalued} 
		{\sc K. Deimling},
		{\it Multivalued differential equations},
		De Gruyter, Berlin, New York, ( 1992).
		
		\bibitem{deng1993exponential}
		{\sc K. Deng},
		{\it Exponential decay of solutions of semilinear parabolic equations with nonlocal initial conditions},
		Journal of Mathematical analysis and applications, {\bf 169}, (2) (1993), 630--637.
		
		\bibitem{du2011initialized}
		{\sc M. Du and Z. Wang},
		{\it Initialized fractional differential equations with Riemann-Liouville fractional-order derivative},
		The European Physical Journal Special Topics, {\bf 193}, (1) (2011), 49--60.
		
		
		\bibitem{heymans2006physical} 
		{\sc N. Heymans and I. Podlubny}
		{\it Physical interpretation of initial conditions for fractional differential equations with Riemann-Liouville fractional derivatives},
		Rheologica Acta, {\bf 45}, (2006), 765--771.
		
		\bibitem{jiang2023topological} 
		{\sc Y. Jiang, Z. Wei,  G. Tang and I. Moroz},
		{\it Topological properties of solution sets for nonlinear evolution hemivariational inequalities and applications},
		Nonlinear Analysis: Real World Applications, {\bf 71}, (2023), 103798.
		
		\bibitem{kalman1963mathematical} 
		{\sc  R. E. Kalman},
		{\it Mathematical description of linear dynamical systems},
		Journal of the Society for Industrial and Applied Mathematics, Series A: Control, {\bf 1}, (2) (1963), 152--192.
		
		\bibitem{kamenskii2011condensing} 
		{\sc  M. I. Kamenskii,  V. V. Obukhovskii and  P. Zecca},
		{\it 	Condensing Multivalued Maps and Semilinear Differential Inclusions in Banach Spaces},
		de Gruyter,Berlin, New York, (2001).
		
		\bibitem{kavitha2021results}
		{\sc K. Kavitha, V. Vijayakumar, A. Shukla, K. S. Nisar and R. Udhayakumar},
		{\it Results on approximate controllability of Sobolev-type fractional neutral differential inclusions of Clarke subdifferential type},
		Chaos, Solitons \& Fractals, {\it 151}, (2021), 111264.
		
		\bibitem{kilbas2006theory} 
		{\sc A. A. Kilbas , H. M. Srivastava, and J. J. Trujillo },
		{\it Theory and applications of fractional differential equations},
		elsevier, (2006).
		
		\bibitem{liang2022existence} 
		{\sc Y. Liang},
		{\it Existence and Approximate Controllability of Mild Solutions for Fractional Evolution Systems of Sobolev-Type},
		Fractal and Fractional, {\bf 6}, (2) (2022), 56.
		
		\bibitem{liu2013approximate} 
		{\sc Z. Liu and M. Bin},
		{\it Approximate controllability for impulsive Riemann-Liouville fractional differential inclusions},
		Abstract and Applied Analysis, {\bf 2013}, (2013), Art. ID 639492.
		
		\bibitem{liu2013controllability} 
		{\sc Z. Liu and X. Li},
		{\it On the controllability of impulsive fractional evolution inclusions in Banach spaces},
		Journal of Optimization Theory and Applications, {\bf 156}, (1) (2013), 167--182.
		
		\bibitem{liu2015approximate} 
		{\sc Z. Liu and X. Li},
		{\it Approximate Controllability of Fractional Evolution Systems with Riemann--Liouville Fractional Derivatives},
		SIAM Journal on Control and Optimization, {\bf 53}, (4) (2015), 1920--1933.
		
		\bibitem{ma2023hilfer} 
		{\sc Y. Ma, C. Dineshkumar, V. Vijayakumar, R. Udhayakumar, A. Shukla and K. S.Nisar},
		{\it Hilfer fractional neutral stochastic Sobolev-type evolution hemivariational inequality: Existence and controllability},
		Ain Shams Engineering Journal, {\bf 14}, (9) (2023), 102126.
		
		\bibitem{mahmudov2003approximate} 
		{\sc N. I. Mahmudov},
		{\it Approximate controllability of semilinear deterministic and stochastic evolution equations in abstract spaces},
		SIAM journal on control and optimization, {\bf 42}, (5) (2003), 1604--1622.
		
		\bibitem{mahmudov2008approximate} 
		{\sc N. I. Mahmudov},
		{\it 	Approximate controllability of evolution systems with nonlocal conditions},
		Nonlinear Analysis: Theory, Methods \& Applications, {\bf 68}, (3) (2008), 536--546.
		
		\bibitem{mohan2024results}
		{\sc M. Mohan Raja, V. Vijayakumar, R. Udhayakumar and K. S. Nisar},
		{\it Results on existence and controllability results for fractional evolution inclusions of order 1< r< 2 with Clarke's subdifferential type},
		{Numerical Methods for Partial Differential Equations} \textbf{40}, e22691 (2024).
		
		\bibitem{papageorgiou1985theory} 
		{\sc  N. S. Papageorgiou},
		{\it On the theory of Banach space valued multifunctions. 1. Integration and conditional expectation},
		Journal of multivariate analysis, {\bf 17}, (2) (1985), 185--206.
		
		\bibitem{PANAGIOTOPOULOS1981335}  
		{\sc  P. D. Panagiotopoulos},
		{\it Non-convex superpotentials in the sense of F.H. Clarke and applications},
		Mechanics Research Communications, {\bf 8}, (6) (1981), 335--340.
		
		\bibitem{panagiotopoulos1993hemivariational} 
		{\sc P. D. Panagiotopoulos},
		{\it Hemivariational inequalities},
		Springer, (1993).
		
		\bibitem{panagiotopoulos2012inequality}
		{\sc  P. D. Panagiotopoulos},
		{\it Inequality Problems in Mechanics and Applications: Convex and nonconvex energy functions},
		Springer Science \& Business Media, (2012).
		
		\bibitem{podlubny2001geometric} 
		{\sc I. Podlubny},
		{\it Geometric and physical interpretation of fractional integration and fractional differentiation},
		Fractional Calculus and Applied Analysis, {\bf 5}, (4) (2002), 367--386.
		
		\bibitem{shu2019approximate} 
		{\sc L. Shu, X. Shu and J. Mao},
		{\it Approximate controllability and existence of mild solutions for Riemann-Liouville fractional stochastic evolution equations with nonlocal conditions of order $1< \alpha< 2$},
		Fractional Calculus and Applied Analysis, {\bf 22}, (4) (2019), 1086--1112.
		
		\bibitem{shi2016study} 
		{\sc X. Shu and Y. Shi},
		{\it A study on the mild solution of impulsive fractional evolution equations},
		Applied Mathematics and Computation, {\bf 273}, (2016), 465--476 .
		
		\bibitem{wang2017approximate} 
		{\sc J. Wang, M. Fe{\v{c}}kan and Y. Zhou},
		{\it Approximate controllability of Sobolev type fractional evolution systems with nonlocal conditions},
		Evolution Equations \& Control Theory, {\bf 6}, (3) (2017), 471--486.
		
		\bibitem{zeng2018class} 
		{\sc S. Zeng, Z. Liu and  S. Migorski},
		{\it A class of fractional differential hemivariational inequalities with application to contact problem},
		Zeitschrift f{\"u}r angewandte Mathematik und Physik, {\bf 69}, (2018), 1--23.
		
		\bibitem{zhou2013existence} 
		{\sc Y. Zhou, L. Zhang, X. H. Shen},
		{\it Existence of mild solutions for fractional evolution equations},
		Journal of Integral Equations and Applications, {\bf 25}, (4) (2013), 557--586.
		
		
	\end{thebibliography}
\end{document}